\documentclass{article}

\usepackage{amsmath}
\usepackage{amsthm}
\usepackage{amssymb}
\usepackage{amsfonts}
\usepackage{latexsym}
\usepackage{delarray}
\usepackage[dvips]{graphics}
\usepackage{epsfig}
\usepackage{color}

\newtheorem{thm}{Theorem}[section]
\newtheorem{lemma}[thm]{Lemma}
\newtheorem{prop}[thm]{Proposition}
\newtheorem{cor}[thm]{Corollary}

\theoremstyle{definition}

\theoremstyle{remark}

\newcommand{\Aut}{\operatorname{Aut}}

\title{The distinguishing number of the direct product and wreath product action}
\author{Melody Chan \\ University of Cambridge\\ Cambridge, England \\ \texttt{melody.chan@aya.yale.edu}}

\begin{document}
\maketitle

\begin{abstract}
Let $G$ be a group acting faithfully on a set $X$.  The distinguishing number of the action of $G$ on $X$, denoted $D_G(X)$, is the smallest number of colors such that there exists a coloring of $X$ where no nontrivial group element induces a color-preserving permutation of $X$.  In this paper, we consider the distinguishing number of two important product actions, the wreath product and the direct product.  Given groups $G$ and $H$ acting on sets $X$ and $Y$ respectively, we characterize the distinguishing number of the wreath product $G~\wr_Y H$ in terms of the number of distinguishing colorings of $X$ with respect to $G$ and the distinguishing number of the action of $H$ on $Y$.  We also prove a recursive formula for the distinguishing number of the action of the Cartesian product of two symmetric groups $S_m \times S_n$ on $[m] \times [n]$.
\end{abstract}

\section{Introduction}

Let $G$ be a group acting faithfully on a set $X$.  For $r \in \mathbb{N}$, an {\em $r$-coloring} of $X$ is  a function $c \colon  X \rightarrow \{1, \ldots, r\}$.  A permutation $\pi$ of $X$ {\em preserves} the coloring $c$ if $c(x^{\pi}) = c(x)$ for all $x \in X$.  A coloring is said to be {\em distinguishing} if the only element in $G$ that induces a color-preserving permutation of $X$ is the identity element.  The {\em distinguishing number} of the action of $G$ on $X$, denoted $D_G(X)$, is the smallest $r$ admitting a distinguishing $r$-coloring of $X$ with respect to the action of $G$.  If there does not exist a distinguishing $r$-coloring of $X$ for any finite $r$, we say that $D_G(X) = \infty$.  

Note that we may equivalently view a distinguishing $r$-coloring of $X$ as a partition $\{X_1, \ldots, X_r\}$ of $X$ into disjoint classes with the property that $G$ intersects the permutation group $X_1 !\times \ldots \times X_r !$ trivially.  The distinguishing number is then the smallest number $r$ admitting such a partition, or $\infty$ if no such $r$ exists.

In \cite{a_c}, Albertson and Collins first introduced the distinguishing number as a property of graphs.  More specifically, the distinguishing number of a graph $M$, denoted $D(M)$, is the smallest number of colors admitting a coloring of the vertices such that the only color-preserving automorphism of $M$ is the identity; thus $D(M) = D_{\Aut(M)}(V(M))$.  The distinguishing number of several families of graphs, including trees, hypercubes, and generalized Petersen graphs, has been computed in \cite{b_c}, \cite{mc_cubes}, \cite{cctcheng}, and \cite{potanka}.  In \cite{jt}, Tymoczko generalized the notion of the distinguishing number to group actions on sets and studied the actions of $S_n$.  In \cite{mc_maximum}, we provided upper bounds for the distinguishing numbers admitted by a large class of groups including nilpotent and supersolvable groups.

We would like to better understand the distinguishing number in the generalized context of group actions introduced by Tymoczko.  To this end, we consider the behavior of the distinguishing number with respect to two natural and important group products: the wreath product and the direct product.  Not only are these products and their associated actions of intrinsic interest, they also allow us to relate the distinguishing number of the action of a large group to the distinguishing numbers of the actions of smaller groups.  

In Section \ref{wreathproduct}, we completely characterize the distinguishing number of the action of the wreath product of two groups on the Cartesian product of their sets.  Our result relates the distinguishing number of the wreath product action to the distinguishing number of one group action and the number of distinct distinguishing colorings of the other group action.  As immediate corollaries, we derive an upper bound for the distinguishing number of imprimitive group actions and a lower bound for the distinguishing number of the lexicographic product of two graphs.  

In Section \ref{directproduct}, we give a recursive formula for the distinguishing number of the direct product of two symmetric groups acting on the direct product of their sets.  This gives an upper bound for the general direct product action.

Our definition of the distinguishing number of a group action differs from the one given in \cite{jt} in that we require the action to be faithful.  This apparent restriction does not actually limit the question being considered.  Given a nonfaithful action of $G$ on $X$, we may consider instead the faithful action of the quotient group $G/Stab(X)$ on $X$, where $Stab(X)$ denotes the elements of $G$ that fix each $x \in X$.  Also, in contrast to both \cite{a_c} and \cite{jt}, we do not require our groups and sets to be finite, simply because there seems to be no reason to do so.  We only note that if $G$ is an infinite group acting faithfully on a set $X$, then $X$ must be infinite as well.

Throughout the paper, we denote group actions by exponentiation on the right.  Thus, the image of an element $x \in X$ under the action of $g \in G$ is denoted $x^g$, and we have $(x^{g_1})^{g_2} = x^{(g_1 g_2)}$ for all $g_1, g_2 \in G$.  The exponentiation notation has the advantage of being relatively intelligble in more complex actions such as the wreath product action.  Also, if $n$ is a positive integer, we use $[n]$ to denote the set $\{1,\ldots,n\}$.

\section{The wreath product action} \label{wreathproduct}

Our main goal in this section is to compute the distinguishing number of the action of the wreath product of two permutation groups on the Cartesian product of the sets upon which they act.

Before defining the wreath product action, let us first recall the definition of the semidirect product of two groups.  Let $A$ and $B$ be our groups, and suppose we have a homomorphism $\phi\colon  B \rightarrow \Aut(A)$.  This homomorphism determines an action of $B$ on $A$ which we will denote by right exponentiation, thus $\phi(b)\colon  a \mapsto a^b$.  Then the semidirect product of $A$ and $B$ according to this action is denoted $A \rtimes_\phi B$ and is the group whose elements are $A \times B$ and whose law of composition is given by
$$
(a_1, b_1) (a_2, b_2) = (a_1 ({a_2}^{b_1^{-1}}), b_1 b_2).
$$
Note that the semidirect product of two groups is not in general uniquely defined, but rather is dependent upon the choice of $\phi$.

Now, let $G$ and $H$ be groups acting faithfully on sets $X$ and $Y$ respectively.  Let $G^Y$ denote the set of functions from $Y$ to $G$.  We equip $G^Y$ with group structure in the following way: given two functions $f_1$ and $f_2$ in $G^Y$, let $f_1 f_2$ be the function given by $f_1 f_2\colon  y \mapsto f_1(y) f_2(y)$.  Note that the identity element of $G^Y$ is the constant 1 function $\mathbf{1}\colon  y \mapsto 1$.  Now, define a homomorphism $\phi\colon  H \rightarrow \Aut(G^Y)$ as follows: for each $h \in H$, we let $(f^h)(y) = f(y^{h^{-1}})$, where $f^h$ denotes the image of $f$ under $h$ according to the right action of $H$ on $G^Y$ determined by $\phi$.  Then the {\em wreath product} of $G$ and $H$ is denoted $G~\wr_Y H$ and is equal to the semidirect product $G^Y \rtimes_\phi H$.  We note that the identity element of this group is $(\mathbf{1}, 1)$.  Finally, we define a right action of $G~\wr_Y H$ on the set $X \times Y$, defined as follows: for any $(x,y) \in X \times Y$ and $(f,h) \in G ~\wr_Y H$ we let $(x,y)^{(f,h)} = (x^{f(y)}, y^h)$.  This action is clearly faithful.  

The wreath product action arises naturally in several important instances.  In order to motivate the ensuing discussion on the distinguishing number of this action, we state a few of them below.

Recall first that an action of a group $K$ on a set $\Omega$ is {\em transitive} if for every $\omega_1, \omega_2 \in \Omega$, there exists $k \in K$ such that ${\omega_1}^k = \omega_2$.  
An equivalence relation $\sim$ on $\Omega$ is {\em $K$-invariant} if $\omega_1 \sim \omega_2$ implies ${\omega_1}^k \sim {\omega_2}^k$ for all $\omega_1, \omega_2 \in \Omega$ and $k \in K$.  Thus we always have two $K$-invariant relations: the universal relation and the relation of equality.  We will call these trivial relations.  A transitive action of $K$ on $\Omega$ is {\em imprimitive} if it admits a nontrivial $K$-invariant relation.  A {\em block of imprimitivity} is an equivalence class under such a relation.

The following proposition tells us that every faithful group action that is transitive but imprimitive is embeddable in a wreath product action.

\begin{prop} \cite[Theorem 2.7]{cameron} \label{p:imprimitive}
Let $K$ be a group acting faithfully, transitively, and imprimitively on $\Omega$.  Let $X$ be a block of imprimitivity of this action, and let $Y = \{ X^k~|~ k \in K \}$ be the set of images of $X$ under the action of $K$.  Let $G$ be the permutation group arising by restricting the setwise stabilizer of $X$ to $X$, and let $H$ be the permutation group that $K$ induces on $Y$.  Then there exists an embedding of the action of $K$ on $\Omega$ into the action of $G~\wr_{Y} H$ on $X \times Y$.
\end{prop}

The significance of Proposition \ref{p:imprimitive} in the context of distinguishing numbers becomes apparent once we state the following lemma.

\begin{lemma} \label{l:subgroup}
Suppose $G$ acts faithfully on $X$.  Let $H$ be a subgroup of $G$ and consider the action of $H$ on $X$ obtained by restricting the action of $G$.  Then $D_H(X) \le D_G(X)$.
\end{lemma}

\begin{proof} If $D_G(X) = \infty$, there is nothing to prove.  Otherwise, there exists a $D_G(X)$-coloring of $X$ such that no nonidentity element of $G$ is color-preserving.  In particular, no nonidentity element of $H \le G$ is color preserving.  So $D_G(X)$ colors suffice to distinguish the action of $H$ on $X$.
\end{proof}

Thus, Proposition \ref{p:imprimitive} tells us that for a faithful, transitive, and imprimitive action of $K$ on $\Omega$ and $G$, $H$, $X$, and $Y$ as defined above,
$$D_K(\Omega) \le D_{G~\wr_{Y} H}(X \times Y).$$
In other words, the distinguishing number of a wreath product action gives an upper bound for the distinguishing number of an imprimitive action embedded in it.  We refer the reader to \cite{cameron} for a more detailed discussion of imprimitive actions and the wreath product.  

As a second example, consider the \emph{wreath product} of two graphs, also called the \emph{graph lexicographic product} or \emph{graph composition}.  Given graphs $\Gamma_1 = (V_1, E_1)$ and $\Gamma_2 = (V_2, E_2)$, the wreath product $\Gamma_1[\Gamma_2]$ is defined to be the graph on vertex set $V_1 \times V_2$ in which two vertices $(v_1, v_2)$ and $(w_1, w_2)$ are connected by an edge if and only if $(v_1, w_1) \in E_1$ or $(v_1 = w_1 \textrm{ and } (v_2, w_2) \in E_2).$
Note that $\Aut(\Gamma_1)~\wr_{V(\Gamma_2)} \Aut(\Gamma_2) \le \Aut(\Gamma_1[\Gamma_2])$, so by Lemma \ref{l:subgroup}, the distinguishing number of the wreath product action gives a lower bound $$D(\Gamma_1[\Gamma_2]) \ge D_{\Aut(\Gamma_1)~\wr_{V(\Gamma_2)} \Aut(\Gamma_2)}(V_1 \times V_2).$$
In \cite{sabidussi}, Sabidussi gives necessary and sufficient conditions for $\Aut(\Gamma_1)~\wr_{V(\Gamma_2)} \Aut(\Gamma_2) = \Aut(\Gamma_1[\Gamma_2])$, in which case this lower bound becomes equality.  His work is generalized in \cite{hemminger} and extended to color digraphs in \cite{d_m}.

With these examples in mind, we now present the main theorem of this section characterizing the distinguishing number of the wreath product action.

\begin{thm}  Let $G$ and $H$ act faithfully on sets $X$ and $Y$ respectively.  For each positive integer $r$, let $n_r$ be the number of distinct distinguishing $r$-colorings of $X$, and let $D_H(Y) = d < \infty$.  Let $S$ be the set $\bigl \lbrace r~|~n_r \ge d|G| \bigr \rbrace$.  Then 
$$
D_{G ~\wr_Y H}(X \times Y) = 
\begin{cases}
\min(S) & \textrm{if~} S \ne \emptyset \\
\infty & \textrm{if~} S = \emptyset
\end{cases}
$$
\end{thm}

\begin{proof}
We will consider the finite and infinite cases separately.

{\it Case:} $S \ne \emptyset$.
Let $k = \min(S)$.  We begin by constructing a distinguishing $k$-coloring of $X \times Y$.  Let $A$ be the set of distinguishing $k$-colorings of $X$, thus $|A| = n_k$.  Now consider the action of $G$ on $A$ defined as $(a^g)(x) = a(x^{g^{-1}})$ for each $a \in A$ and $g \in G$.  Each $a \in A$ is distinguishing, so it has trivial stabilizer and orbit length $|G|$.  Therefore the number of orbits of the action of $G$ on $A$ is $|A|/|G| = n_k/|G|$.  Since $n_k/|G| \ge d$, we may choose $d$ distinguishing $k$-colorings of $X$ in pairwise disjoint orbits.  Call these $k$-colorings $a_1, \ldots, a_d$.  Now, let $b$ be a distinguishing $d$-coloring of $Y$ with respect to the action of $H$ (the existence of which follows from the assumption that $D_H(Y) = d$).  Let $C\colon  X \times Y \rightarrow \{1, \ldots, k\}$ be given by $C\colon  (x,y) \mapsto a_{b(y)}(x)$.  We claim that $C$ is a distinguishing $k$-coloring of $X \times Y$.

Suppose that $(f,h) \in G~\wr_Y H$ preserves $C$.  We wish to show that $(f,h) = (\mathbf{1}, 1)$.  First, we show that $h$ preserves the coloring $b$.  We know that for each $(x,y) \in X \times Y$, $C(x,y) = C((x,y)^{(f,h)}) = C(x^{f(y)}, y^h)$.  By definition of $C$, we have $a_{b(y^h)}(x^{f(y)}) = a_{b(y)}(x)$, and so ${a_{b(y^h)}}^{f(y)^{-1}} = a_{b(y)}$.  Therefore ${a_{b(y^h)}}$ and $a_{b(y)} \in A$ are in the same orbit under the action of $G$.  But we chose $a_i$ and $a_j$ to be in different orbits if $i \ne j$.  Therefore, $b(y^h) = b(y)$ for each $y \in Y$.  Thus, $h$ permutes the elements of $Y$ in a way that preserves the coloring $b$.  Since $b$ is a distinguishing coloring of $Y$, we must have $h = 1$.

Now we show that $f(y) = 1$ for each $y \in Y$.  Indeed, the fact that $(f,1)$ preserves $C$ gives that $a_{b(y)}(x^{f(y)}) = a_{b(y)}(x)$ for each $x \in X$, $y \in Y$.  Thus $f(y)$, considered as a permutation of $X$, preserves the coloring $a_{b(y)}$ of $X$.  Since each $a_{b(y)}$ is a distinguishing coloring of $X$, we have $f(y) = 1$ for each $y$.  Thus, $f = \mathbf{1}$, and $(f,h) = (\mathbf{1}, 1)$, and we conclude that $C$ is a distinguishing $k$-coloring of $X \times Y$.

It remains to be shown that every distinguishing coloring of $X \times Y$ uses at least $k$ colors.  Suppose that $C'$ is a distinguishing $l$-coloring of $X \times Y$; we will show that $l \ge k$.  For each $y \in Y$, let $a_y\colon  X \rightarrow \{1, \ldots, l\}$ be given by $a_y\colon  x \mapsto C'(x,y)$ for each $x \in X$.  Now, we claim each $a_y$ must be a distinguishing $l$-coloring of $X$.  For if $g \in G$ preserves $a_y$, let $f \in G^Y$ be given by $f(y') = g$ if $y' = y$ and $f(y') = 1$ otherwise.  Then $(f,1) \in G~\wr_Y H$ preserves $C'$, and since $C'$ is a distinguishing coloring of $X \times Y$, we have $f = \mathbf{1}$ and $g = 1$.  So $a_y$ is a distinguishing $l$-coloring of $X$ for each $y \in Y$.  

Now let $A'$ be the set of distinguishing $l$-colorings of $X$, thus $|A'| = n_l$.  Let $\mathbf{O} = \{ O \textrm{ an orbit of the action of } G \textrm{ on } A'~|~a_y \in O \textrm{ for some } y \in Y\}$.  Let $|\mathbf{O}| = d'$, and write $\mathbf{O} = \{O_1, \ldots, O_{d'} \}$.  Now let $b\colon  Y \rightarrow \{1, \ldots, d'\}$ be given by $b\colon  y \mapsto i$ if $a_y$ is in orbit $O_i$.  The function $b$ is well-defined since the orbits are disjoint.  We claim that $b$ is a distinguishing $d'$-coloring of $Y$.  To verify this claim, suppose $h \in H$ preserves $b$.  This means that $b(y^h) = b(y)$ for each $y$, so $a_{y^h}$ and $a_y$ are in the same orbit of $G$ on $A'$ for each $y$.  Then for each $y \in Y$, let $g_y \in G$ be the element taking $a_{y^h}$ to $a_y$, thus ${(a_{y^h})}^{g_y} = a_y$.  Now let $f\colon  Y \rightarrow G$ be given by $f\colon  y \mapsto {g_y}^{-1}$.  We claim that $(f,h) \in G~\wr_Y H$ preserves the coloring $C'$.  Indeed, for each $(x,y) \in X \times Y$, we have $C'((x,y)^{(f,h)}) = C'(x^{{g_y}^{-1}}, y^h) = a_{y^h}(x^{{g_y}^{-1}}) = {(a_{y^h})}^{g_y}(x) = a_y(x) = C'(x,y)$.  Since $C'$ was assumed to be a distinguishing coloring of $X \times Y$, we have $(f,h) = (\mathbf{1}, 1)$ so $h = 1$.  Thus, $b$ is a distinguishing $d'$-coloring of $Y$.  Since $D_H(Y) = d$, we have $d' \ge d$.

Finally, we note again that since each $a' \in A'$ is a distinguishing coloring, it has trivial stabilizer and orbit length $|G|$ under the natural action of $G$.  Then the number of orbits of this action is $|A'|/|G| = n_l/|G|$.  But the number of orbits is at least $|\mathbf{O}| = d' \ge d$, so $n_l \ge d'|G| \ge d|G|$.  Since $k$ was the minimum number such that $d|G| \le n_k$, we have $l \ge k$.  Thus every distinguishing coloring of $X \times Y$ must use at least $k$ colors, and we conclude that $D_{G ~\wr_Y H}(X \times Y) = k$. 

{\it Case:} $S = \emptyset$.  We have seen above that if $X \times Y$ had a distinguishing $l$-coloring for finite $l$, then $n_l \ge d|G|$.  Since no such $n_l$ exists, we must have $D_{G ~\wr_Y H}(X \times Y) = \infty$.

\end{proof}

In \cite{cctcheng}, Cheng shows that $n_r$, the number of distinct distinguishing $r$-colorings of $X$ with respect to the action of $G$, is always a monic polynomial in $r$ of degree $|X|$.  She furthermore gives a recursive formula that computes $n_r$ in the case that $G$ is the automorphism group of a tree acting on vertex set $X$.  It seems that $n_r$ is in general difficult to compute.  However, we give an explicit formula for the special cases $G = S_n$ and $G = A_n$ in the following corollaries.

\begin{cor} \label{c:sn_wreath} Suppose $H$ acts faithfully on $Y$ with distinguishing number $d$.  Then $D_{S_n ~\wr_Y H}([n] \times Y) = \min \{ r~|~\binom{r}{n} \ge d \}$.
\end{cor}

\begin{proof} The distinguishing number of $S_n$ on $[n]$ is clearly $n$.  So for a fixed $r$, there are $\binom{r}{n}n! = \binom{r}{n}|S_n|$ distinguishing $r$-colorings of $[n]$.
\end{proof}

Before presenting the next corollary, we first consider the action of the alternating group $A_n$ on $[n]$.
\begin{lemma} $D_{A_n}([n]) = n-1$.  
\end{lemma}
\begin{proof}
Given $n-1$ distinct colors, we may color the elements of $[n]$ such that only $1$ and $2$ in $[n]$ share a color.  Since the transposition $(1~2)$ is not in $A_n$, no nontrivial permutation preserves this coloring.  On the other hand, with $n-2$ or fewer colors available, either at least three elements $a$, $b$, and $c \in [n]$ share a color, in which case $(a~b~c)$ is color-preserving, or we have at least two pairs of elements, $a$ and $b$ colored identically and $c$ and $d$ colored identically, in which case $(a~b)(c~d)$ is color-preserving.  Thus $A_n$ has distinguishing number $n-1$ in its natural action.
\end{proof}

\begin{cor} Suppose $H$ acts faithfully on $Y$ with distinguishing number $d$.  Then $D_{A_n ~\wr_Y H}([n] \times Y) = \min \{ r~|~(n-1) \binom{r}{n-1} + 2 \binom{r}{n} \ge d \}$.
\end{cor}

\begin{proof}
There are $\binom{r}{n-1} \binom{n}{2} (n-1)!$ ways to choose a distinguishing coloring of $[n]$ using $n-1$ distinct colors.  There are $\binom{r}{n} n!$ ways to choose a distinguishing coloring of $[n]$ using $n$ colors.  Thus, we require $\binom{r}{n-1} \binom{n}{2} (n-1)! + \binom{r}{n} n! \ge d |A_n| = d (\frac{n!}{2})$, whence the result follows.

\end{proof}

\section{The direct product action} \label{directproduct}

Given groups $G$ and $H$ acting faithfully on sets $X$ and $Y$ respectively, what is the distinguishing number of $G \times H$ acting on $X \times Y$?  In this section, we answer the above question when $G$ and $H$ are the full symmetric groups $S_m$ and $S_n$ in their natural actions on $[m]$ and $[n]$.  
Because the distinguishing number of the natural action of the symmetric group is easily computable, and the direct product is such a simple way to combine two group actions, it is quite surprising that the characterization of $D_{S_m \times S_n}([m] \times [n])$ is so complex.
Yet this action is of particular interest because it gives an upper bound for all finite direct product actions.  Specifically, if $G$ and $H$ are groups acting faithfully on sets $X$ and $Y$ of cardinalities $m$ and $n$ respectively, then by Lemma \ref{l:subgroup}, we have
$$D_{G \times H}(X \times Y) \le D_{S_m \times S_n}([m] \times [n]).$$
We will see that $D_{S_m \times S_n}([m] \times [n])$ is often very small, so the upper bound given above is often a useful one.  In general, $D_{S_m \times S_n}([m] \times [n])$ does not depend so much on the absolute size of $m$ and $n$ as it depends on their relative size.  The farther apart $m$ and $n$ are, the greater the distinguishing number of the corresponding action is.

We first prove the well-definedness of a function that will be used in the main theorem.

\begin{lemma} \label{l:welldefined}
Fix $k \ge 2$.  Then there is a unique function $f_k\colon \{2, 3, \ldots \} \rightarrow \mathbb{N}$ satisfying:

(1) if $m \le k$ then $f_k(m) = 1$, and

(2) if $m > k$ then $f_k(m)$ is the smallest integer $t$ such that $1 < t < m$ and $m \le k^t - f_k(t)$.
\end{lemma}

\begin{proof}

We fix $k \ge 2$ and proceed by induction on $m$ to show that $f_k(m)$ is well-defined.  If $m \le k$, then $f_k(m) = 1$.  Now suppose $m > k$ and assume inductively that $f_k(i)$ is well-defined for $1 < i < m$.  It suffices to show that the set 
$$S_{k,m} = \{t~|~1 < t < m \textrm{ and } m \le k^t - f_k(t)\}$$ 
is nonempty.  Note that if $t < m$, then $f_k(t) \le t-1$ by the inductive hypothesis.  Thus, $k^{m-1} - f_k(m-1) \ge k^{m-1} - (m-2)$.  Furthermore, one may check that $k^{m-1} - (m-2) \ge m$ for each $k \ge 2$ and $m \ge 2$.  So 
$$m \le k^{m-1} - f_k(m-1).$$  Thus, $m-1 \in S_{k,m}$ and therefore $f_k(m) = \min (S_{k,m})$ is well-defined for each $m$.
\end{proof}

\begin{thm} \label{t:fn}  
Fix $m \ge 2$ and $n \ge 1$ and let $f_k(m)$ be defined as in Lemma \ref{l:welldefined}.  Then the set
$$T_{m,n} = \{ k \ge 2~|~f_k(m) \le n \le k^m - f_k(m)\}$$
is nonempty, and
$$D_{S_m \times S_n}([m] \times [n]) = \min(T_{m,n}).$$
\end{thm}

Note that we restrict $m \ge 2$ only for convenience in the proof; if $m = 1$ then the action of $S_m \times S_n$ is isomorphic to the action of $S_n$ on $[n]$ and has distinguishing number $n$.  Also, it is interesting that the symmetry between $m$ and $n$ is not at all obvious from the formulation of Theorem \ref{t:fn}.

Throughout the proof, we will regard the set $[m] \times [n]$ as a grid of $m$ rows and $n$ columns.  An element of $S_m \times S_n$ acts on this grid as a permutation of the rows followed by a permutation of the columns.  We can think of a coloring of the $m \times n$ grid as composed of $n$ column colorings.  In a distinguishing coloring, each of the $n$ column colorings must be distinct, for otherwise two identically colored columns could be transposed to produce a nontrivial color-preserving permutation.  We also note that if every column in a grid has a distinct coloring, the only color-preserving group element that leaves the rows unchanged is the identity element.

We will show that $f_k(m)$ gives the smallest number $n$ such that the $m \times n$ grid has a distinguishing $k$-coloring.  Furthermore, we will prove that the $m \times x$ grid has a distinguishing $k$-coloring precisely when $x$ is between $f_k(m)$ and $k^m - f_k(m)$.  The proof of this fact will proceed by induction on $k$ with base case $k = 2$.  The theorem then follows.

We first show that if $n$ is too large with respect to $m$ and $k$, then $k$ colors do not suffice for a distinguishing coloring.

\begin{lemma} \label{l:2}
If $n \ge k^m$ then the $m \times n$ grid does not have a distinguishing $k$-coloring.
\end{lemma}
\begin{proof}

Let $c$ be a $k$-coloring of the $m \times n$ grid.  Then each column must have a distinct coloring.  There are $k^m$ possible column colorings, so we must have $n = k^m$ and each column coloring is used exactly once.  But then for any nontrivial row permutation $\sigma$, there exists a nontrivial column permutation $\tau$ such that $\sigma$ and $\tau$ induce identical changes in the coloring $c$.  Then $(\sigma, \tau^{-1})$ is a nontrivial color-preserving permutation of the $m \times n$ grid, contradicting the assumption that $c$ is distinguishing.

\end{proof}

\begin{lemma} \label{l:3}
Suppose $1 \le n \le k^m - 1$.  Then the $m \times n$ grid has a distinguishing $k$-coloring if and only if the $m \times (k^m - n)$  grid has a distinguishing $k$-coloring.
\end{lemma}
\begin{proof}
Proving one direction suffices by symmetry.  Let $c$ be a distinguishing $k$-coloring of the $m \times n$ grid.  Then $c$ uses exactly $n$ of the $k^m$ possible column colorings.  Let $c'$ be a coloring of the $m \times (k^m - n)$ grid where each of the remaining $k^m - n$ column colorings is used exactly once.   We claim that $c'$ is distinguishing.  

Any nontrivial row permutation $\sigma$ applied to the coloring $c$ of the $m \times n$ grid must introduce some column coloring not occurring in $c$ (and therefore occurring in $c'$).  For otherwise, $\sigma$ would only have permuted the column colorings of $c$ and so some column permutation $\tau$ could restore $c$, contradicting the assumption that $c$ is distinguishing.  But then $\sigma$ applied to the coloring $c'$ of the $m \times (k^m - n)$ grid must introduce some column coloring not in $c'$.  This shows that $\sigma$ cannot be the row component of a permutation that preserves $c'$.  So only permutations that leave the rows unchanged could possibly preserve $c'$.  But we have already noted that only the identity element falls into this category.  This proves that $c'$ is distinguishing.

\end{proof}

The next two lemmas give some conditions under which $k$ colors do not suffice for a distinguishing coloring.

\begin{lemma} \label{l:4}
If $1 \le n < f_k(m)$ then the $m \times n$ grid does not have a distinguishing k-coloring.
\end{lemma}
\begin{proof}
Fix $k$ and proceed by induction on $m$, with base cases $m \le k$ that are vacuously true since $f_k(m) = 1$ in this case.  Also note that if $n = 1$ then our assumption on $n$ gives that $f_k(m) > 1$, so $m > k$ and there does not exist a distinguishing $k$-coloring of the $m \times 1$ grid.  So we may assume $n \ge 2$.  Now, $n < f_k(m)$ implies that $m > k^n - f_k(n)$ by definition of $f_k(m)$.  So $k^n - m < f_k(n)$.  If $m \ge k^n$ then the $m \times n$ grid does not have a distinguishing $k$-coloring by Lemma \ref{l:2}.  So suppose instead that $1 \le k^n - m < f_k(n)$.  Since $n < f_k(m) < m$, we may apply the inductive hypothesis to conclude that there does not exist a distinguishing $k$-coloring of the $n \times (k^n -m)$ grid.  Then by Lemma \ref{l:3}, there does not exist a distinguishing $k$-coloring of the $n \times m$ grid, and therefore of the $m \times n$ grid.
\end{proof}

\begin{lemma} \label{l:5}
If $k^m - f_k(m) < n \le k^m - 1$ then the $m \times n$ grid does not have a distinguishing $k$-coloring.
\end{lemma}
\begin{proof} This result follows from Lemma \ref{l:3} and Lemma \ref{l:4}.
\end{proof}

The next lemma is the key result that allows us to construct distinguishing colorings of large grids from distinguishing colorings of small ones.

\begin{lemma} \label{l:55}

Suppose $c$ is a distinguishing $k$-coloring of the $m \times n$ grid and $N$ is the number of column colorings that cannot be obtained via a row permutation from some column coloring in $c$.  Then for every $l$ satisfying $n \le l \le n+N$, the $m \times l$ grid also has a distinguishing $k$-coloring.

\end{lemma}

\begin{proof}

We construct a distinguishing $k$-coloring of the $m \times l$ grid, where $n \le l \le n+N$, as follows.  Color the leftmost $n$ columns as they are colored in $c$; at most $N$ columns remain.  Color these remaining columns with distinct column colorings, none of which can be obtained from some column coloring in $c$ via a row permutation.  Call the resulting coloring $c'$.  We claim that this is a distinguishing coloring of the $m \times l$ grid.  Any nontrivial row permutation $\sigma$ must take some column coloring in $c$ to one not in $c$; otherwise some column permutation $\tau$ could restore $c$, contradicting that $c$ is distinguishing.  This means that $\sigma$ must also take some column coloring in $c'$ to one not in $c'$, for none of the additional column colorings in $c'$ can be obtained via a row permutation.  This shows that no nontrivial row permutation can be part of a color-preserving permutation of $c'$.  Since $c'$ gives a distinct coloring for each column, it must therefore be distinguishing.

\end{proof}

The next three lemmas give some conditions that guarantee the existence of a 2-coloring.  This case will provide the base case of a proof that proceeds by induction on the number of colors.

First, we note that a coloring is distinguishing if each column coloring is distinct and each row contains a different number of color 1 entries.

\begin{lemma} \label{l:distinctrows}

Let $c$ be a 2-coloring of a grid such that each column has a different coloring and each row has a different number of color 1 entries.  Then $c$ is distinguishing.

\end{lemma}

\begin{proof}

A transformation that preserves $c$ cannot permute the rows nontrivially, since each row has a different number of color 1 entries.  So it cannot permute the columns nontrivially either, since each column has a distinct coloring.

\end{proof}

\begin{lemma} \label{l:square}
The $m \times m$ grid has a distinguishing 2-coloring.
\end{lemma}

\begin{proof}
Let 
$$c((i,j)) = \begin{cases} 1 & \textrm{ if } i < j, \\ 2 & \textrm{ otherwise.} \end{cases}$$
Then Lemma \ref{l:distinctrows} gives that $c$ is distinguishing.
\end{proof}

\begin{lemma} \label{l:case2}
For each $m \ge 2$, if $f_2(m) \le n \le 2^m - f_2(m)$ then there exists a distinguishing 2-coloring of the $m \times n$ grid.
\end{lemma}
\begin{proof}
We proceed by induction on $m$.  If $m = 2$ then we note that the $2\times 1$, $2 \times 2$, and $2 \times 3$ grids have distinguishing $2$-colorings 
$$\left( \begin{array}{c} 1\\2 \end{array} \right), \left( \begin{array}{c c} 1 & 1 \\1 & 2 \end{array} \right), \textrm{ and } \left( \begin{array}{c c c} 1&1&2\\1&2&2 \end{array} \right)$$ 
respectively.  
Now, suppose that for each $2 \le i < m$, we know that the $i \times n$ grid has a distinguishing 2-coloring if $f_2(i) \le n \le 2^i - f_2(i)$.  We wish to show that this property holds for $m$.  Note that it suffices to prove that $m \times n$ grid has a distinguishing 2-coloring if $f_2(m) \le n \le (2^m)/2 = 2^{m-1}$, for the remaining case $2^{m-1} < n \le 2^m - f_2(m)$ must then hold by the symmetry provided by Lemma \ref{l:3}.  We will make repeated use of this condition.

{\em Case 1}: $f_2(m) \le n < m$.  Now, $m$ is at least 3 so $n \ge f_2(m) \ge 2$.  Applying the inductive hypothesis for $n$, it suffices to prove that $f_2(n) \le m \le 2^n - f_2(n)$.  The first inequality is certainly true since $f_2(n) < n < m$.  As for the second, note that by its definition, $f_2(m)$ satisfies $m \le 2^{f_2(m)} - f_2(f_2(m))$.  Now, one may show inductively that $f_2(x)$ increases by at most 1 when $x$ increases by 1.  Then $2^x - f_2(x)$ is an increasing function of $x$, so since $n \ge f_2(m)$, we have $m \le 2^n -f_2(n)$ as desired.

{\em Case 2}: $m \le n \le 2^{m-1}$.  We break the analysis into further cases.

{\em Case 2.1}: $m = 3$.  We need to check the cases $n = 3$ and $n = 4$.  If $n=3$ then the $m \times n$ grid has a distinguishing 2-coloring by Lemma \ref{l:square}.  If $n=4$ then the coloring 
$$\left( \begin{array}{c c c c} 1&1&1&2 \\ 1&1&2&2 \\ 1&2&2&2 \end{array} \right)$$
is a distinguishing 2-coloring of the $3 \times 4$ grid by Lemma \ref{l:distinctrows}.

{\em Case 2.2}: $m = 4$.  The coloring
$$\left( \begin{array}{c c c c} 1&1&1&2\\1&2&2&1\\2&1&2&2\\2&2&2&2 \end{array}\right)$$
gives a distinguishing 2-coloring of the $4 \times 4$ grid by Lemma \ref{l:distinctrows}.  Each column contains either 1 or 2 entries of color 1, so there are $2^4 - (\binom{4}{1} + \binom{4}{2}) = 6$ possible column colorings that cannot be obtained from any of the above 4 column colorings via a row permutation.  Lemma \ref{l:55} then tells us that the $4 \times n$ grid has a distinguishing 2-coloring if $4 < n \le 10$, which is more than sufficient since $2^{4-1} = 8$.

{\em Case 2.3}: $m = 5$.  The coloring 
$$\left( \begin{array}{c c c c} 1&1&1&1\\1&1&1&2\\1&2&2&1\\2&1&2&2\\2&2&2&2 \end{array}\right)$$ 
gives a distinguishing 2-coloring of the $5 \times 4$ grid by Lemma \ref{l:distinctrows}.  Each column contains either 2 or 3 entries of color 1, so Lemma \ref{l:55} tells us that the $5 \times n$ grid has a distinguishing 2-coloring if $4 < n \le 4 + 2^5 - (\binom{5}{2} + \binom{5}{3}) = 16$.

{\em Case 2.4}: $m \ge 6$ is even.  First, we note that the $m \times m$ grid has a distinguishing 2-coloring by Lemma \ref{l:square}.  We construct a distinguishing 2-coloring $c$ of the $m \times (m+1)$ grid as follows: let $$c((i,j)) = \begin{cases}
2 & \textrm{ if } (j = 1) \textrm{ and } (i = 2 \textrm{ or } i \ge m/2 + 2), \\
2 & \textrm{ if } (j = 2) \textrm{ and } (3 \le i \le m/2 + 1 \textrm{ or } i = m), \\
2 & \textrm{ if } (3 \le j \le m/2+1) \textrm{ and } (j \le i \le m/2 \textrm{ or } i \ge m - j + 2), \\
2 & \textrm{ if } (m/2 +1 < j \le m+1) \textrm{ and } (i = m+2 - j \textrm{ or } (i > m/2 \textrm{ and } i \ne 3m/2 + 2 - j)), \\
1 & \textrm{ otherwise.}
\end{cases}
$$
As an example, $c$ is shown below for $m = 10$.
$$\left( \begin{array}{c c c c c c c c c c c c} 
1&1&1&1&1&1&1&1&1&1&2\\2&1&1&1&1&1&1&1&1&2&1\\1&2&2&1&1&1&1&1&2&1&1\\1&2&2&2&1&1&1&2&1&1&1\\1&2&2&2&2&1&2&1&1&1&1\\1&2&1&1&1&2&2&2&2&2&1\\
2&1&1&1&2&2&2&2&2&1&2\\2&1&1&2&2&2&2&2&1&2&2\\2&1&2&2&2&2&2&1&2&2&2\\2&2&2&2&2&2&1&2&2&2&2
\end{array}\right)$$
The coloring $c$ has the property that each column has a distinct coloring and each row contains a different number of entries of color 1.  (In fact, the $i^{th}$ row contains $m+1 - i$ entries of color 1).  So by Lemma \ref{l:distinctrows}, $c$ is distinguishing.  One may check that each column has $m/2$ entries of color 1 and $m/2$ entries of color 2, so by applying Lemma \ref{l:55} to $c$, we find that the $m \times n$ grid has a distinguishing 2-coloring if $m < n \le (m+1) + 2^m - \binom{m}{m/2}$.  
One may check that $\binom{m}{m/2} \le 2^{m-1}$, so that $(m+1) + 2^m - \binom{m}{m/2} \ge 2^{m-1}$ as desired.

{\em Case 2.5}: $m > 6$ is odd.  Then $m-1 \ge 6$ is even, so let $c$ be the distinguishing 2-coloring of the $(m-1) \times m$ grid as given above in Case 2.4.  Let $c'$ be a 2-coloring of the $m \times m$ grid obtained by adding a row of entries colored 1 to the top of $c$.  Then in $c'$, each column has a distinct coloring, and each row contains a different number of entries of color 1.  (In fact, the $i^{th}$ row contains $m+1 - i$ entries of color 1).  So by Lemma \ref{l:distinctrows}, $c'$ is distinguishing.  Furthermore, each row column contains $(m+1)/2$ entries of color 1 and $(m-1)/2$ entries of color 2, so by Lemma \ref{l:55}, the $m \times n$ grid has a distinguishing 2-coloring if $m < n \le m + 2^m - \binom{m}{(m-1)/2}$.  One may check that $\binom{m}{(m-1)/2} \le 2^{m-1}$, so that $m + 2^m - \binom{m}{(m-1)/2} \ge 2^{m-1}$ as desired.  This completes the proof of Lemma \ref{l:case2}.
\end{proof}

Lemma \ref{l:case2} will serve as a base case in the following induction on the number of colors.
\begin{lemma} \label{l:casen}
For each $k \ge 2$ and $m \ge 2$, if $f_k(m) \le n \le k^m - f_k(m)$, then the $m \times n$ grid has a distinguishing $k$-coloring.
\end{lemma}
\begin{proof}
We proceed by induction on $k$.  The case $k = 2$ is precisely Lemma \ref{l:case2}.  Now fix $k > 2$.  Our inductive hypothesis will be that for each $m \ge 2$, the $m \times n$ grid has a distinguishing $(k-1)$-coloring if $f_{k-1}(m) \le n \le (k-1)^m - f_{k-1}(m)$.  We wish to prove that for each $m \ge 2$, the $m \times n$ grid has a distinguishing $k$-coloring if $f_k(m) \le n \le k^m - f_k(m)$.

We first claim that it is sufficient to prove that for each $m \ge 2$, the $m \times n$ grid has a distinguishing $k$-coloring if $f_k(m) \le n < f_{k-1}(m)$.  We check the other cases below.  If $f_{k-1}(m) \le n \le (k-1)^m - f_{k-1}(m)$ then the $m \times n$ grid has a distinguishing $(k-1)$-coloring by the inductive hypothesis, which is certainly also a distinguishing $k$-coloring.  Now consider in particular a distinguishing $(k-1)$-coloring $c$ of the $m \times ((k-1)^m - f_{k-1}(m))$ grid, which we may view as a $k$-coloring where the color $k$ is never used.  There are $k^m - (k-1)^m$ column colorings that use color $k$ at least once and hence cannot be obtained via a permutation from any column coloring in $c$.  Then by Lemma \ref{l:55}, the $m \times n$ grid has a distinguishing $k$-coloring if 
$$(k-1)^m - f_{k-1}(m) \le n$$
and
\begin{eqnarray*}
n &\le& (k-1)^m - f_{k-1}(m) + k^m - (k-1)^m \\
&=& k^m - f_{k-1}(m).
\end{eqnarray*}
Next, we note that by Lemma \ref{l:3}, the $m \times l$ grid has a distinguishing $k$-coloring for all $l$ such that $k^m - f_{k-1}(m) < l \le k^m - f_k(m)$ if and only if the $m \times n$ grid has a distinguishing $k$-coloring for all $n$ such that $f_k(m) \le n < f_{k-1}(m)$.  Thus we need only consider the case $f_k(m) \le n < f_{k-1}(m)$.

To prove our claim that for each $m \ge 2$, the $m \times n$ grid has a distinguishing $k$-coloring if $f_k(m) \le n < f_{k-1}(m)$, we proceed again by induction, this time on $m$, with base case $2 \le m \le k$.  If $2 \le m < k$ (and $k > 2$ is still fixed), then $f_k(m) = f_{k-1}(m) = 1$ and the condition $f_k(m) \le n < f_{k-1}(m)$ is vacuous.  So the statement is (vacuously) true for $m < k$.  If $m = k$, then $f_k(m) = 1 \le n < f_{k-1}(m) = 2$ so $n = 1$ and there does indeed exist a distinguishing $k$-coloring of the $m \times 1$ grid.
  
Now assume inductively that each $i$ with $2 \le i < m$ has the property that the $i \times n$ grid has a distinguishing $k$-coloring if $f_k(i) \le n < f_{k-1}(i)$.  We wish to show that the $m \times n$ grid has a distinguishing $k$-coloring if $f_k(m) \le n < f_{k-1}(m)$.
Only the case $m > k$ remains to be considered.  

If $m > k$, then $f_k(m) \ge 2$ so $n \ge f_k(m) \ge 2$.  So since $n < f_{k-1}(m) < m$, it suffices to prove that $f_k(n) \le m \le k^n - f_k(n)$, for the inductive hypothesis then gives that the $n \times m$ grid has a distinguishing $k$-coloring.  The first inequality is certainly true since $f_k(n) < n < m$.  As for the second, note that by its definition, $f_k(m)$ satisfies $m \le k^{f_k(m)} - f_k(f_k(m))$.  Now, one may show that $f_k(x)$ increases by at most 1 when $x$ increases by 1.  Then $k^x - f_k(x)$ is an increasing function of $x$, so since $n \ge f_k(m)$, we have $m \le k^n -f_k(n)$ as desired.  We conclude that if $f_k(m) \le n \le k^m - f_k(m)$, then we have a distinguishing $k$-coloring of the $m \times n$ grid.

\end{proof}

We combine these results below to prove Theorem \ref{t:fn}.
\begin{proof}
Fix $k, m \ge 2$.  
If $n$ satisfies $f_k(m) \le n \le k^m - f_k(m)$ then the $m \times n$ grid has a distinguishing $k$-coloring by Lemma \ref{l:casen}.  On the other hand, if $n < f_k(m)$ or $n > k^m - f_k(m)$ then the $m \times n$ grid does not have a distinguishing $k$-coloring by Lemmas \ref{l:2}, \ref{l:4}, and \ref{l:5}.
So there exists a distinguishing $k$-coloring of the $m \times n$ grid if and only if $f_k(m) \le n \le k^m - f_k(m)$.  Then by definition of the distinguishing number, $D_{S_m \times S_n}([m] \times [n]) = \min \{ k \ge 2 ~|~ f_k(m) \le n \le k^m - f_k(m)\}$.  

Note that we only needed to consider colorings using at least 2 colors because $m \ge 2$ implies that $S_m \times S_n$ is nontrivial and so acts with distinguishing number at least 2.  
\end{proof}

For $k$ fixed, the function $f_k(m)$ grows approximately logarithmically with $m$.  Thus, the expression $k^m - f_k(m)$ is dominated by $k^m$ for large $m$.  So for a fixed $m$ sufficiently large, the distinguishing number $D_{S_m \times S_n}([m] \times [n])$ grows approximately like the function $\root m \of n$ when $n$ becomes large.  Table \ref{to10} gives $D_{S_m \times S_n}([m] \times [n])$ for $m$ and $n$ between 1 and 10.

\begin{table}[h] 
\begin{tabular}{c|c c c c c c c c c c}
$m,n$ & 1 & 2 & 3 & 4 & 5 & 6 & 7 & 8 & 9 & 10 \\
\hline
1   & 1 & 2 & 3 & 4 & 5 & 6 & 7 & 8 & 9 & 10 \\
2   & 2 & 2 & 2 & 3 & 3 & 3 & 3 & 3 & 4 & 4 \\
3   & 3 & 2 & 2 & 2 & 2 & 2 & 3 & 3 & 3 & 3 \\
4   & 4 & 3 & 2 & 2 & 2 & 2 & 2 & 2 & 2 & 2 \\
5   & 5 & 3 & 2 & 2 & 2 & 2 & 2 & 2 & 2 & 2 \\
6   & 6 & 3 & 2 & 2 & 2 & 2 & 2 & 2 & 2 & 2 \\
7   & 7 & 3 & 3 & 2 & 2 & 2 & 2 & 2 & 2 & 2 \\
8   & 8 & 3 & 3 & 2 & 2 & 2 & 2 & 2 & 2 & 2 \\
9   & 9 & 4 & 3 & 2 & 2 & 2 & 2 & 2 & 2 & 2 \\
10  & 10& 4 & 3 & 2 & 2 & 2 & 2 & 2 & 2 & 2
\end{tabular}
\caption{The distinguishing number of the action of $S_m \times S_n$ on $[m] \times [n]$.}
\label{to10}
\end{table}

\section{Discussion and open questions}

Section \ref{wreathproduct} shows that we can characterize the distinguishing number of $G~\wr_Y H$ on $X \times Y$ if we have information about $n_r$, the number of distinct distinguishing $r$-colorings of the action of $G$ on $X$, for each $r$.  We saw that $n_r$ could be computed when $G$ was the automorphism group of a tree or when $G = A_n$ or $S_n$.  It would be useful to find other examples of group actions for which the $n_r$ can be explicitly computed.  Regarding Section \ref{directproduct}, we ask whether one can provide a closed formula for the distinguishing number of $S_m \times S_n$ in its action on $[m] \times [n]$.  We also ask what the distinguishing number of the general direct product action is.  

There are many interesting questions to ask regarding the distinguishing number of group actions.  In \cite{mc_maximum}, for example, we define $\overline{D}(G)$ to be the maximum distinguishing number admitted by a given group $G$.  Given two groups $G$ and $H$ such that $H \le G$, we ask whether it must be the case that $\overline{D}(H) \le \overline{D}(G)$.  We also ask for a characterization of the set $$\{D_G([n])~|~G \textrm{ is a transitive subgroup of } S_n\}.$$  Note that we require our group $G$ to be transitive, for otherwise each distinguishing number $k$ between $1$ and $n$ could be achieved by taking a subgroup of $S_n$ that fixes each $k+1, k+2, \ldots, n$ and whose action on $1, \ldots, k$ is isomorphic to the action of $S_k$.

We refer the reader to \cite{mc_maximum} for other open questions.

\section{Acknowledgments}

This research was conducted at the University of Minnesota Duluth Research Experience for Undergraduates, while the author was a student at Yale University, and supported by the National Science Foundation (DMS-0137611).  The author would like to thank Melanie Wood and Philip Matchett for many ideas and suggestions on drafts of this paper and to Joseph Gallian for his support and encouragement.

\end{document}